\documentclass[12pt]{amsart}

\usepackage{amsmath, graphicx, cite, enumerate, comment}
\usepackage{amssymb, amsthm, amscd}
\usepackage{latexsym}
\usepackage{amssymb}
\usepackage[english]{babel}
\usepackage{amsmath,amssymb}
\usepackage{amsfonts}
\usepackage{amscd}
\usepackage{cite}
\usepackage{marvosym}
\usepackage{mathrsfs}
\usepackage{amsthm}
\usepackage{hyperref}
\usepackage{verbatim}

\usepackage{tikz} \usetikzlibrary{arrows,shapes,patterns} % AC: You need this command for everything below to work
\usepackage{lipsum} % for sample

\newtheorem{thm}{Theorem}

\newtheorem{prop}[thm]{Proposition}

\theoremstyle{remark}
\newtheorem{rmk}[thm]{Remark}

\theoremstyle{definition}

\numberwithin{equation}{section}

\newcommand{\R}{\mathbb{R}}

\newcommand{\be}{\begin{equation}}
\newcommand{\ee}{\end{equation}}
\newcommand{\bd}{\begin{displaymath}}
\newcommand{\ed}{\end{displaymath}}

       \thanks{\textsl{Mathematics Subject Classification (MSC 2010):} Primary 53A10; Secondary 53C42, 49Q05.}
     \title[A note on the index of closed minimal hypersurfaces of flat tori]{A note on the index of closed minimal hypersurfaces of flat tori}
     \author{Lucas Ambrozio, Alessandro Carlotto and Ben Sharp}
     \address{ \noindent L. Ambrozio: Imperial College London, South Kensington Campus, London SW7 2AZ, United Kingdom, \textit{E-mail address: l.ambrozio@imperial.ac.uk} \newline \newline \indent A. Carlotto: ETH - Department of Mathematics, R\"amistrasse 101, 8092 Z\"urich,Switzerland, \textit{E-mail address: alessandro.carlotto@math.ethz.ch} \newline \newline 
     \indent B. Sharp: University of Warwick, Gibbet Hill Rd, Coventry CV4 7AL, United Kingdom \textit{E-mail address: b.sharp@warwick.ac.uk}}
     
\begin{document}
     	
     	\begin{abstract} Generalizing earlier work by Ros in ambient dimension three, we prove an affine lower bound for the Morse index of closed minimal hypersurfaces inside a flat torus in terms of their first Betti number (with purely dimensional coefficients).
     			\end{abstract}
     	
     	\maketitle       
							
	\section{Introduction}\label{sec:intro}
	
	Motivated by a variety of recent constructions of closed minimal hypersurfaces in positively curved Riemannian manifolds, and by the associated natural classification questions, we presented in \cite{AmbCarSha} a study of the relation between their Morse index and their Betti numbers, namely those pieces of data respectively encoding the most basic analytic and topological information concerning the hypersurface in question. This relation is in fact the object of a conjecture due to Schoen and Marques-Neves \cite{Mar,Nev} that can be stated as follows: in any closed manifold of positive Ricci curvature there is a linear lower bound of the Morse index of a minimal hypersurface $M^n$ in terms of its first Betti number, that is to say
	\begin{equation}\label{eq:conj}
	index(M)\geq C b_1(M)
	\end{equation}
	for some constant $C$ only depending on the ambient manifold. We refer the reader to the introduction of \cite{AmbCarSha} for a broader contextualization of this problem and for a discussion of the various cases for which we could verify this conjecture. \\
	\indent On the other hand, it is straightforward to observe that inequality \eqref{eq:conj} cannot hold true in the special but fundamental case of flat manifolds, as is seen by considering totally geodesic $n$-dimensional tori inside $(n+1)$-dimensional flat tori (in which case one has $index(M)=0$ and $b_1(M)=n$ for any $n\geq 2$).
	In such a setting, the best one can hope for is instead an affine bound with a negative additive constant on the right-hand side. Up to now, an estimate of that sort was only obtained for $n=2$ by A. Ros \cite{Ros}. The scope of this note is to prove the following generalization of such a result:
		
			\begin{thm} \label{thm:main} 
				Let $M^n$ be a closed minimal hypersurface in a $(n+1)$-dimensional flat torus. Assume there is a point $p$ in $M^n$ where all principal curvatures are distinct. Then
				\begin{equation*}
				index(M) \geq \frac{2}{n(n+1)} (b_{1}(M) - (2n-1)).
				\end{equation*}
				\indent If $n=2$ or $n=3$, then the above inequality holds true without the assumption on the principal curvatures. 
			\end{thm}
			
			Clearly, there is a natural correspondence between minimal hypersurfaces in flat tori and $(n+1)$-periodic minimal hypersurfaces in the Euclidean space $\mathbb{R}^{n+1}$, a topic that has been thoroughly investigated with several interesting results: we refer the reader to the article by W. Meeks \cite{Mee} and to the survey by D. Hoffman \cite{Hoff} as well as references therein for further details about the classical case $n=2$. Remarkably, lots of interesting examples in $\mathbb{R}^3$ are actually known, the most basic ones being provided by the Schwarz $P$ and $D$ surfaces, the latter class ensuring that the estimate of Theorem \ref{thm:main} is actually sharp since in that case one can find, in a fundamental domain, that the Morse index equals one (this follows from the work of M. Ross \cite{Ross}) while the topology is that of a genus three orientable surface. By contrast, the construction of periodic minimal hypersurfaces in ambient dimension at least four
			is a fascinating theme of current research: see in particular the recent work by Choe and Hoppe \cite{ChHop} for certain higher-dimensional analogues of the aforementioned classical surfaces and related comments therein. \\

				\indent In applying Theorem \ref{thm:main} to obtain classification results, it is also useful to know that non-trivial minimal hypersurfaces inside a flat torus must have sufficiently large first Betti number:

	\begin{thm}(Cf. Theorem 1 in \cite{Kel}, and Theorem 4.1 in \cite{ChoFra})\label{thm:toptor}
		Let $M^n$ be a closed minimal hypersurface in a $(n+1)$-dimensional flat torus. Then $b_{1}(M) \geq n+1$ unless $M^n$ is a flat totally geodesic $n$-dimensional torus.
	\end{thm}
	
	\indent This fact follows from a more general statement that goes back at least to E. Kelly \cite{Kel}, but see also Theorem 1 in \cite{NagSmy} for an interesting generalization to harmonic maps and Theorem 4.1 in \cite{ChoFra} for a broad analysis of topological restrictions imposed by the existence of minimal immersions into manifolds of Ricci curvature bounded from below. For the sake of completeness, a simple and direct proof of Theorem \ref{thm:toptor} is presented in Subsection \ref{subs:proofmains}. \\

	\indent In \cite{Ros}, Ros proved that a non-orientable, compact stable minimal surface immersed in a flat three-torus $T^3$ has the topology of a Klein bottle with a handle (cf. Theorem 7 therein). Analogously, as a simple combined application of the two theorems above one can prove that if $M^3\subset T^4$ is stable (but not totally geodesic) then either $b_1=4$ or $b_1=5$ and $M^3$ is non-orientable as can be checked directly by means of the second variation formula. The question of classifying all such stable minimal hypersurfaces remains a challenging open problem. \\

	%\indent At a technical level, our main result follows from an upper bound on the number of harmonic forms $\omega$ on $M^n$ such that all ambient components of the 2-form $N^{\flat}\wedge \omega$
	%(that is to say: all functions of the form $\langle N^{\flat}\wedge \omega, \theta_i\wedge \theta_j\rangle$ for $\left\{\theta_1,\ldots, \theta_n\right\}$ a global parallel frame for the ambient manifold) are Jacobi fields. In turn, this relies on the one hand on an interesting identity characterizing the value of the Jacobi operator applied to any such function (Proposition 3) and on the other on the analysis of the rigidity case (Proposition 4) that is inspired by a recent article of Chao Li \cite{Li}, which we improve and adapt to our setting. 
	
		\textsl{Acknowledgments}. The authors wish to express their gratitude to Andr\'e Neves for his interest in this work and for a number of enlightening conversations.
		During the preparation of this article,	L.A. was supported by prof. Neves' ERC Start Grant PSC and LMCF 278940.

	\section{Proofs}\label{sec:proofs}
	
	\subsection{Notations and ancillary results} 
	
	\
	
	\
	
	Throughout this note, we consider the ambient manifold $T^{n+1}:=\R^{n+1}/\Gamma$, where $\Gamma$ is a lattice group of maximal rank, endowed with its (flat) Riemannian metric $\langle \cdot,\cdot\rangle$ and the associated Levi-Civita connection $D$. Furthermore, we let $M^n$ denote a closed, embedded minimal hypersurface in $T^{n+1}$, that is to say a smooth, closed hypersurface of vanishing mean curvature. For the sake of simplicity, and in order to streamline our arguments, we shall assume that $M^n$ is orientable or, equivalently, two-sided and we let $N$ denote a choice of its unit normal vector field. The case when $M^n$ is one-sided, for which the statement of Theorem \ref{thm:main} still holds, is discussed in Remark \ref{rem:onesid}. The induced Levi-Civita connection on the submanifold $M^n$ will always be denoted by $\nabla$, while $\Delta$ stands for the Laplace-Beltrami operator. Lastly, we convene that the second fundamental form of $M^n\subset T^{n+1}$ is defined by the formula $A(X,Y)=\langle D_X Y, N\rangle$ for any pair of smooth vector fields $X, Y$ along $M^n$. \\

	\indent In our setting the Morse index can be defined as follows: the sections of the normal bundle of $M^n\subset T^{n+1}$ can be identified with the set of smooth functions $\phi$ on $M^n$, and the second variation of the area functional is given by the quadratic form
	\begin{equation*} \label{defindexfrom}
	Q(\phi,\phi) = \int_{M} (|\nabla \phi|^2 -|A|^2\phi^2) dM,
	\end{equation*}
	\noindent so that the number of \textsl{negative variations} is encoded in the spectrum of the Jacobi operator 
	\begin{equation*} \label{defijacobi}
	J_M \phi = \Delta\phi +|A|^2\phi
	\end{equation*} 
 and the index of $M^n$ is by definition the number of (strictly) negative eigenvalues of $J_M$. \\

	\indent In a closed Riemannian manifold $(M^n,g)$, the Ho\-dge-\-La\-pla\-ce operator is the second order differential operator $\Delta_{p}$ acting on $p$-forms defined by
	 \begin{equation*}
	 \Delta_p = dd^{*} + d^{*}d
	 \end{equation*}  
	 \noindent where $d : \Omega^{p}(M) \rightarrow \Omega^{p+1}(M)$ is the exterior differential and $d^{*}:  \Omega^{p}(M) \rightarrow \Omega^{p-1}(M) $ is the formal adjoint of $d$, defined with respect to the metric $g$. A $p$-form $\omega$ is called harmonic when $\Delta_{p}\omega = 0$ and we let $\mathcal{H}^p(M,g)$ denote the vector space of harmonic $p$-forms on $(M^n,g)$. When $M^n$ is closed, $\omega$ is harmonic if and only if it is closed and co-closed, that is to say when both $d\omega=0$ and $d^{*}\omega = 0$ hold true. Hodge's Theorem asserts that in a closed Riemannian manifold one has the isomorphism $\mathcal{H}^1(M,g)\simeq H^1(M;\mathbb{R})$ so that the dimension of the space of harmonic $1$-forms coincides with the first Betti number of the manifold. 
	 Also, we will use the (special) Bochner-Weitzenb\"ock formula relating the Hodge-Laplace operator with the usual (rough) Laplacian on $1$-forms:
	 \begin{equation} \label{eq:bochner}
	 \Delta_{1}\omega = - \Delta\omega + Ric^{M}(\omega^{\sharp},\cdot).
	 \end{equation}

	 	\indent In this note we employ the usual musical isomorphisms to pass from vectors to $1$-forms, see Remark 2.1 in \cite{AmbCarSha} for further details.\\

\indent	The proof of our main result relies on the following proposition:
	
	\begin{prop}\label{prop:ident} Let $M^n$ be a closed, orientable, minimal hypersurface in $T^{n+1}$ and let $\omega$ be a harmonic $1$-form on $M^n$. For every parallel $2$-form $\theta$ on $T^{n+1}$ one has the identity
		\begin{equation} \label{eqmain}
		\Delta\langle N^{\flat}\wedge\omega,\theta\rangle  + |A|^2\langle N^{\flat}\wedge\omega,\theta\rangle  = -2\sum_{i,j=1}^{n}A(E_i,E_j)\langle E^{\flat}_{j}\wedge \nabla_{E_i}\omega,\theta\rangle,
		\end{equation}
		\noindent where the expression on the right-hand side is globally defined as it does not depend on the particular choice of local orthonormal frame $\{E_i\}$ on $M^n$.
	\end{prop}
	
	\begin{rmk}
		The reader may want to compare this assertion with Lemma 1 in \cite{Ros}, where similar computations were performed for the coordinates of $\omega$ rather than $N^{\flat}\wedge\omega$ (the two choices being essentially equivalent only if $n=2$). 
	\end{rmk}
	
	\begin{proof}
		Let $\{E_i\}$ be a local orthonormal frame on $M^n$, which is geodesic at a point $p$ in $M^n$ (that is to say $(\nabla_{E_i}E_j)(p)=0$ for all $i,j$). We have
		\begin{equation*}
		D_{E_i}(N^{\flat}\wedge\omega) = D_{E_i}N^{\flat}\wedge \omega + N^{\flat}\wedge D_{E_i}\omega = -\sum_{j=1}^nA(E_i,E_j)E^{\flat}_{j}\wedge \omega + N^{\flat}\wedge\nabla_{E_i}\omega.
		\end{equation*}
		\indent Hence, exploiting the fact that $\theta$ is parallel, one has at the point $p$
		\begin{align*}
		\Delta\langle N^{\flat}\wedge\omega,\theta\rangle  = {} & \sum_{i=1}^n E_i E_i \langle N^{\flat}\wedge\omega,\theta\rangle  = \sum_{i=1}^n E_i \langle D_{E_i}(N^{\flat}\wedge\omega),\theta\rangle  \\
		= {} & \sum_{i=1}^n E_i \langle -\sum_{j=1}^n A(E_i,E_j)E^{\flat}_{j}\wedge \omega + N^{\flat}\wedge\nabla_{E_i}\omega,\theta\rangle  \\
		= {} & -\sum_{i,j=1}^n E_i (A(E_i,E_j)) \langle E^{\flat}_{j}\wedge \omega,\theta\rangle - \sum_{i,j=1}^n A(E_i,E_j)\langle D_{E_i}E^{\flat}_j\wedge\omega,\theta\rangle  \\
		& - \sum_{i,j=1}^n A(E_i,E_j)\langle E^{\flat}_j\wedge D_{E_i}\omega,\theta\rangle + \sum_{i=1}^n\langle D_{E_i}(N^{\flat}\wedge\nabla_{E_i}\omega), \theta\rangle.
		\end{align*}
		\indent Since $M^n$ is minimal, by the Codazzi equation for flat ambient manifolds  (evaluating, once again, at $p$)
		\begin{equation*}
		\sum_{j=1}^n E_{j}(A(E_i,E_j))=\sum_{j=1}^n E_i (A(E_j,E_j))= E_i H = 0.
		\end{equation*}
		\indent Thus, we can deduce 
		\begin{align*}
		\Delta\langle N^{\flat}\wedge\omega,\theta\rangle  = & - \sum_{i,j=1}^n A(E_i,E_j)A(E_i,E_j)\langle N^{\flat}\wedge\omega,\theta\rangle - \sum_{i,j=1}^n A(E_i,E_j)\langle E^{\flat}_j\wedge \nabla_{E_i}\omega,\theta\rangle  \\
		& - \sum_{i,j=1}^n A(E_i,E_j)A(E_i,\omega^{\sharp})\langle E^{\flat}_j\wedge N^{\flat},\theta\rangle - \sum_{i,j=1}^n A(E_{i},E_{j})\langle E^{\flat}_j\wedge\nabla_{E_i}\omega, \theta\rangle  \\
		& + \sum_{i=1}^n \langle N^{\flat}\wedge\nabla_{E_i}\nabla_{E_i}\omega,\theta\rangle.
		\end{align*}
		\indent We can rewrite the above as
		\begin{equation*}
		\Delta\langle N^{\flat}\wedge\omega,\theta\rangle  + |A|^2\langle N^{\flat}\wedge\omega,\theta\rangle = \langle \Delta\omega,i_N\theta\rangle  + \langle  A\circ A(\omega^{\sharp}),(i_N\theta)^{\sharp}\rangle- 2\sum_{i,j=1}^n A(E_{i},E_{j})\langle  E^\flat_j\wedge\nabla_{E_i}\omega, \theta\rangle.
		\end{equation*}
		\indent Now, the Gauss equation for a minimal hypersurface in flat ambient manifolds yields
		\begin{equation*}
		Ric^M(\omega^{\sharp},(i_N\theta)^{\sharp}) = -\langle A\circ A(\omega^{\sharp}),(i_N\theta)^{\sharp}\rangle,
		\end{equation*}
		while, on the other hand, the Bochner formula \eqref{eq:bochner} for harmonic $1$-forms on $M^n$ reads
		\begin{equation*}
		\langle \Delta\omega,i_{N}\theta\rangle  = Ric^M(\omega^{\sharp},(i_{N}\theta)^{\sharp})
		\end{equation*}
		so that formula (\ref{eqmain}) follows at once. 
	\end{proof}

	\indent Inspired by Chao Li (see \cite{Li}, Proposition 5.1), we prove the following result.
	
	\begin{prop}\label{prop:vanish}
		Let $M^n$ be a closed minimal hypersurface in a $(n+1)$-dimensional flat torus. Assume there is a point $p$ in $M^n$ where all the principal curvatures are distinct. The set of all harmonic $1$-forms $\omega$ on $M^n$ such that 
		\begin{equation} \label{eq:jacobi}
		\sum_{i,j=1}^{n}A(E_i,E_j) E^{\flat}_j \wedge \nabla_{E_i}\omega = 0
		\end{equation}
		has dimension at most $2n-1$.
	\end{prop} 
	
	Before proceeding with its proof, we want to remind the reader of a general fact about harmonic forms: \textsl{in a Riemannian manifold $(M^n,g)$ a $1$-form $\omega\in\Omega^1(M^n)$ is closed and co-closed if and only if $\nabla\omega$ is a symmetric trace-free tensor.}

	\begin{proof}
		
		By virtue of our assumption on $M^n$, we can certainly find a positive number $\rho$ (smaller than the injectivity radius of $M^n$ at $p$) such that the principal curvatures of $M^n$ are all distinct in the geodesic ball $B_{\rho}(p)$ and, furthermore, there exists a local orthonormal frame $\left\{E_1,\ldots, E_n\right\}$ diagonalizing the second fundamental form $A$ at every point of such ball (namely: $A(E_i,E_j)=k_i\delta_{ij}$ for all $i,j=1,\ldots, n$ with $k_1<k_2<\ldots<k_n$).
		Now, given a harmonic $1$-form $\omega$, equation \eqref{eq:jacobi} and the fact that the tensor $\nabla\omega$ is symmetric imply
		\begin{equation*} 
		(k_{i}-k_{j})(\nabla_{E_i}\omega)(E_j) = 0 \quad \text{for all} \quad i,j=1,\ldots, n,
		\end{equation*}
		so that
		\begin{equation} \label{eq:jacobi2}
		(\nabla_{E_i}\omega)(E_j) = 0 \quad \text{for all} \quad i\neq j
		\end{equation}
		on the whole geodesic ball in question. In particular, since $\nabla\omega$ is also trace-free we deduce that the functions $\nabla\omega(E_1, E_1),\ldots,\nabla\omega (E_{n-1}, E_{n-1})$ completely determine the tensor $\nabla\omega$ on $B_{\rho}(p)$.\\
		\indent Now, let us consider the functions defined on $B_{\rho}(p)$ by
		\begin{equation*}
		\phi_i (q)	= \begin{cases}
		\omega (E_i)(q) &  \textrm{for} \ 1\leq i\leq n \\
		\nabla\omega (E_{i-n},E_{i-n})(q) & \textrm{for} \ n+1\leq i\leq 2n-1.
		\end{cases}
		\end{equation*}
		We claim that if $\omega$ satisfies \eqref{eq:jacobi} then for every $q\in B_{\rho}(p)$ the values $(\phi_1,\ldots, \phi_{2n-1})(q)$ are uniquely determined by $(\phi_1,\ldots, \phi_{2n-1})(p)$, hence the space of harmonic forms in $B_{\rho}(p)$ satisfying \eqref{eq:jacobi} has dimension at most $2n-1$ and the general statement over $M^n$ follows by unique continuation.

			\indent To check the claim, we proceed as follows: given $q\in B_{\rho}(p)$ let $\gamma:[0,\tau]\to M^n$ be the only geodesic connecting $p$ to $q$ in $B_{\rho}(p)$ and consider the functions (of one real variable) obtained by restriction of $(\phi_1,\ldots, \phi_{2n-1})$ along $\gamma$, namely set
		\begin{equation*}
		f_i(t)=\phi_i(\gamma(t)), \ \ \textrm{for} \ i=1,\ldots, 2n-1.
		\end{equation*}
		Then, we claim that $(f_1,...,f_{2n-1})$  solves a linear ODE system (in normal form) so that (by Cauchy-Lipschitz) the value at $p$ uniquely determines the value along the curve $\gamma$, which is enough to check the claim above. \\
		\indent Notice that $\gamma'(t)=\sum_{j=1}^n \alpha^j(t)E_{j}(\gamma(t))$ where the coefficients $\alpha^j, \ j=1,\ldots,n$ are smooth and bounded (since the frame $\left\{E_1,\ldots, E_n\right\}$ is orthonormal) hence by linearity
		\begin{equation}
		\frac{d}{dt}f_i(t)=\nabla_{\gamma'(t)}\phi_i(\gamma(t))=\sum_{j=1}^n \alpha^j(t)\nabla_{E_j}\phi_i (\gamma(t))
		\end{equation}
		so that it suffices for our scopes to check that for every choice of the indices $1\leq i\leq 2n-1$ and $1\leq j\leq n$ the function $\nabla_{E_j}\phi_i$ can be expressed as a linear combination of $\phi_1,\ldots,\phi_{2n-1}$, with smooth coefficients. \\
	    \indent First of all, for $i\leq n$ we have
		\begin{equation}\label{eq:firstcomps}
		\nabla_{E_j}\phi_i=\delta_{ij}\phi_{n+i}+\sum_{k=1}^n \Gamma_{ji}^k \phi_k.
		\end{equation}
		This can be justified as follows:
		\begin{align*}
		\nabla_{E_j}(\omega(E_i))= (\nabla_{E_j}\omega)(E_i)+\omega(\nabla_{E_j}E_i)=\nabla\omega(E_i,E_j)+\sum_{k=1}^n \Gamma_{ji}^k\omega(E_k)= \delta_{ij}\phi_{n+i}+\sum_{k=1}^n \Gamma_{ji}^k\phi_k
			\end{align*}
			where the last steps relies on equation \eqref{eq:jacobi2}.
			Also, observe that for $i=n$ the RHS of \eqref{eq:firstcomps} needs to be suitably interpreted, namely with $-\sum_{l=1}^{n-1}\phi_{n+l}$ in lieu of $\phi_{2n}$.

		On the other hand, the differential equation for $\phi_{n+i}$ takes for $1\leq i\leq n-1$ one of the following three forms:\newline

	\noindent\textsl{Case 1: $1\leq j\leq n-1,  i\neq j$}
		\begin{equation*}\label{eq:secondcomps1}
		\nabla_{E_j}\phi_{n+i}=
		(2\Gamma^i_{ji}-\Gamma^i_{ij})\phi_{n+i}-\Gamma^j_{ii}\phi_{n+j}-\sum_{k=1}^n R_{ikij}\phi_k. 
		 	\end{equation*}
		\textsl{Case 2: $j=n$}
		\begin{equation*}\label{eq:secondcomps2}
\nabla_{E_j}\phi_{n+i}=(2\Gamma_{ni}^i-\Gamma_{in}^i)\phi_{n+i}+\sum_{k=1}^{n-1}\Gamma_{ii}^k \phi_{n+k}-\sum_{k=1}^n R_{ikin}\phi_k.
		\end{equation*}
		\textsl{Case 3: $j=i$}
		\begin{equation*}\label{eq:secondcomps3}
		\nabla_{E_j}\phi_{n+i}= \sum_{k=1, k\neq i}^n \Gamma^i_{kk}\phi_{n+i}-\sum_{k=1, k\neq i}^{n-1}(2\Gamma_{ik}^k-\Gamma_{ki}^k)\phi_{n+k}+ (2\Gamma_{in}^n -\Gamma_{ni}^n)\sum_{k=1}^{n-1}\phi_{n+k}+\sum_{k=1}^n R_{ki}\phi_k.\nonumber
		\end{equation*}
		
		Indeed, for all $j\neq i$ by the Ricci formula for commuting covariant derivatives (no summation on repeated indexes $i$ and $j$) we obtain
			\begin{align*}
 \nabla_{E_j}& (\nabla\omega(E_i,E_i))= \nabla^2\omega(E_i,E_i,E_j) +2\nabla\omega(\nabla_{E_j}E_i,E_i) \\
		= & \nabla^2\omega(E_i,E_j,E_i)-\sum_{k=1}^{n} R_{ikij}\omega(E_{k})+2\Gamma_{ji}^{i}\nabla\omega(E_i,E_i) \\
		= &\nabla_{E_{i}}(\nabla\omega(E_i,E_j))-\nabla\omega(\nabla_{E_i}E_i,E_j)- \nabla\omega(E_i,\nabla_{E_i}E_j) -\sum_{k=1}^{n}R_{ikij}\omega(E_{k})+ 2\Gamma_{ji}^{i}\nabla\omega(E_i,E_i)\\
		= & -\Gamma^j_{ii}\phi_{n+j}+(2\Gamma_{ji}^i-\Gamma^i_{ij})\phi_{n+i}-\sum_{k=1}^n R_{ikij}\phi_k
		\end{align*}
		 where the very last equality was implied by \eqref{eq:jacobi2}. Thereby, the first equation is verified and the second one follows along similar lines exploiting the fact that the tensor $\nabla\omega$ is trace free.
		 Lastly, following the same pattern, one has for all $i=1,\ldots,n-1$ 
		 \begin{align*}
	\nabla_{E_{i}}& (\nabla\omega(E_i,E_i)) =\sum_{k=1, k\neq i}^{n}-E_i(\nabla\omega(E_k,E_k)) 
		  =\sum_{k=1, k\neq i}^{n-1} \Gamma_{kk}^{i}\nabla\omega(E_i,E_i) -(2\Gamma_{ik}^{k}-\Gamma_{ki}^k)\nabla\omega(E_k,E_k)\\
		  & +\Gamma_{nn}^{i}\nabla\omega(E_i,E_i)-(2\Gamma_{in}^{n}-\Gamma_{ni}^n)\nabla\omega(E_n,E_n)+\sum_{k=1, k\neq i}^{n-1} \sum_{l=1}^{n}R_{klki}\omega(E_{l})+\sum_{k=1}^{n}R_{nkni}\omega(E_{k})\\
		  = & \sum_{k=1, k\neq i}^n \Gamma^i_{kk}\phi_{n+i}-\sum_{k=1, k\neq i}^{n-1}(2\Gamma_{ik}^k-\Gamma_{ki}^k)\phi_{n+k}+ (2\Gamma_{in}^n-\Gamma_{ni}^n)\sum_{k=1}^{n-1}\phi_{n+k}+\sum_{k=1}^n R_{ki}\phi_k.
		 \end{align*}
		 \indent This finishes the proof of the claim and thus the whole argument.
		\end{proof}

	\subsection{Proofs of Theorem \ref{thm:main} and Theorem \ref{thm:toptor}}\label{subs:proofmains}
	
	\
	
	\
	
	As stated in the introduction, we first present a short proof of Theorem \ref{thm:toptor}.
	
	\begin{proof}
		Let $\mathcal{V}$ be the space of all parallel $1$-forms on the flat $(n+1)$-dimensional torus. This space consists precisely of the forms $df$, where $f$ is a linear function on the universal cover of the flat torus (viz. $(n+1)$-dimensional Euclidean space). In particular, its dimension is $n+1$. \\
		\indent Since $M^n$ is minimal and the elements of $\mathcal{V}$ are parallel, the restriction of any $df\in \mathcal{V}$ is a harmonic $1$-form. Moreover, $df=0$ on $M^n$ if and only if $M^{n}$ is contained in the quotient of level sets of the linear map $f$ by the action of the lattice subgroup of translations of the Euclidean space that generates the flat torus in question. Our assertion follows at once. 
	\end{proof}
	
	We then deduce from Proposition \ref{prop:ident} and Proposition \ref{prop:vanish} our main result, Theorem \ref{thm:main}. To that scope, we need to recall a rigidity theorem proved by Do Carmo and Dajczer \cite{doCDaj} for minimal hypersurfaces $M^n$ in $\mathbb{R}^{n+1}$ such that some principal curvature has multiplicity at least $n-1$ at all points: any such hypersurface must be part of a higher-dimensional catenoid, or flat (see Corollary 4.4 in their paper). In particular, in a four dimensional torus ($n=3$), the assumption on the principal curvatures given in the statement of Theorem \ref{thm:main} (namely the assumption that all principal curvatures are pairwise distinct) holds for all closed minimal hypersurfaces that are not totally geodesic.
	
	\begin{proof} Given an orthonormal basis $\{\theta_1,\ldots,\theta_{n+1}\}$ of parallel $1$-forms on $T^{n+1}$, let $M^n$ be a closed, orientable minimal hypersurface of a flat $(n+1)$-dimensional torus whose principal curvatures are all distinct at least at one point. Let us denote by $k$ its Morse index, and by $\{\phi_q\}_{q=1}^{k}$ an $L^2$-orthonormal basis of eigenfunctions of the Jacobi operator $J_{M}=\Delta+|A|^2$ of $M^n$ generating the eigenspace where the operator is negative definite. Let then $\Phi$ denote the linear map defined by
		\begin{equation*}
		\begin{matrix}
		\Phi : & \mathcal{H}^{1}(M^n) & \rightarrow  & \mathbb{R}^{n(n+1)k/2} \\
		&  \omega & \mapsto & \left[\int_{M} \langle N^{\flat}\wedge \omega, \theta_{i}\wedge \theta_{j} \rangle\phi_{q} d\mu\right],
		\end{matrix}
		\end{equation*}
		where $1\leq i < j \leq n+1$ and $q$ varies from $1$ to $k$. Clearly,
		\begin{equation*}
		dim  \mathcal{H}^{1}(M,g) \leq dim Ker(\Phi) + \frac{n(n+1)}{2}k.
		\end{equation*}
		 \indent Since $\mathcal{H}^{1}(M,g)\simeq H^{1}(M;\mathbb{R})$, Theorem \ref{thm:main} will follow once we analyse the dimension of the kernel of the map $\Phi$ and show that indeed
		$dim Ker(\Phi)\leq 2n-1$.\\
		\indent Let $\omega$ be an element of the kernel of the map $\Phi$. This precisely means that every function $u_{ij}=\langle N^{\flat}\wedge \omega, \theta_{i}\wedge\theta_{j} \rangle$ is $L^2$-orthogonal to each of the first $k$ eigenfunctions of $J_M$. Since $index(M)=k$, we must have 
		\begin{equation*}
		Q(u_{ij},u_{ij})\geq \lambda_{k+1}\int_{M}u_{ij}^2d\mu \geq 0 \quad \text{for all} \quad 1\leq i<j \leq n+1,
		\end{equation*}
		by the standard variational characterization of eigenvalues. Hence, thanks to Proposition \ref{prop:ident} we have
		 \begin{align*}
		 0 \leq & \sum_{1\leq i< j \leq n+1}Q(u_{ij},u_{ij}) = -\sum_{1\leq i< j \leq n+1}\int_{M}u_{ij}J_{M}(u_{ij})d\mu \\ 
		 & =2\int_M\sum_{1\leq i< j \leq n+1}\langle N^{\flat}\wedge \omega, \theta_i\wedge \theta_j\rangle\langle \sum_{k,l=1}^{n}A(E_k,E_l)E^{\flat}_l\wedge\nabla_{E_k}\omega,\theta_i\wedge\theta_j\rangle\,d\mu \\
		 & =2\int_{M}\sum_{k,l=1}^n A(E_k,E_l)\langle E^{\flat}_l\wedge\nabla_{E_k}\omega,N^{\flat}\wedge\omega\rangle d\mu=0,
		 \end{align*}
		 the last equality relying on the fact that trivially, by orthogonality, $i_{N}(E^{\flat}_l\wedge\nabla_{E_k}\omega)=0$ for any choice of the indices $k$ and $l$. 
		 It follows that $\langle N^{\flat}\wedge\omega,\theta_i\wedge\theta_j\rangle$ are all eigenfunctions of the Jacobi operator $J_M$, associated to the eigenvalue $\lambda_{k+1} = 0$. By Proposition \ref{prop:ident}, this implies that $\omega$ satisfies equation \eqref{eq:jacobi}. Thus, to complete the proof, it is enough to invoke Proposition \ref{prop:vanish}, which ensures that the dimension of $Ker(\Phi)$ cannot exceed $2n-1$. \\
		 \indent Lastly, to obtain an unconditional result when $n=2,3$, we need to observe that when the condition on the principal curvatures of $M^n$ is not fulfilled, then $M^n$ is totally geodesic and therefore a stable $n$-dimensional torus, in which case our inequality is also satisfied. This claim is clear when $n=2$ and is a consequence of the result by do Carmo and Dajczer when $n=3$.
		  Thereby the proof is complete.
			\end{proof}

	\begin{rmk}\label{rem:onesid}
	Let us discuss the modifications needed to obtain this index bound for non-orientable minimal hypersurfaces inside $T^{n+1}$. First of all, by orientability of $T^{n+1}$ we know that any such $M^{n}$ is also one-sided. In this case, it is customary to introduce the two-sheeted covering $\pi: \hat{M}^n\to M^n$ (given by couples $(x,N)$ where $x\in M^n$ and $N$ attains one of the two possible choices for the unit normal of $M^n$ at $x$) and the associated two-sided immersion $\iota:\hat{M}^n\to T^{n+1}$ (with a well-defined unit normal field $\hat{N}$ given by $\hat{N}(x,N)=N$). As discussed (for instance) in Subsection 2.3 of \cite{AmbCarSha} one can then consider the restriction of the Jacobi operator of $\hat{M}^n$ to the space of \textsl{odd} functions, namely to all $u:\hat{M}^n\to\mathbb{R}$ satisfying $u\circ\tau = -u$ for $\tau:\hat{M}^n\to\hat{M}^n$ the natural deck transformation of the covering in question. Hence, the Morse index of $M^n$ is defined to be the number of negative eigenvalues of $J_{\hat{M}}$ acting on odd functions, in the sense above. \\
	\indent These comments being made, our arguments go through in the non-orientable case as well without any substantial modification once it is checked that for any harmonic form $\omega\in\mathcal{H}^1(M,g)$ one has that the test functions $u_{ij}=\langle \hat{N}^{\flat}\wedge\pi^{\ast}(\omega), \iota^*(\theta_i)\wedge\iota^*(\theta_j) \rangle$ are indeed odd, which is done in the proof of Theorem A of \cite{AmbCarSha}.
	\end{rmk}

	\bibliographystyle{amsbook}

\end{document}